\documentclass[pdflatex,sn-mathphys-num]{sn-jnl}
\usepackage[T1]{fontenc}
\usepackage{amsmath,amssymb,amsfonts,amsthm}
\usepackage{microtype}
\usepackage{needspace}
\usepackage{xcolor}
\newcommand{\pyo}[1]{\mathrel{\perp_{P,#1}}}
\theoremstyle{thmstyleone}
\newtheorem{theorem}{Theorem}[section]
\newtheorem{lemma}[theorem]{Lemma}
\newtheorem{proposition}[theorem]{Proposition}
\newtheorem{corollary}[theorem]{Corollary}
\theoremstyle{thmstyletwo}

\newcommand{\R}{\mathbb R}
\newcommand{\N}{\mathbb N}
\newcommand{\spn}{\operatorname{span}}
\newcommand{\norm}[1]{\left\lVert#1\right\rVert}
\newcommand{\mut}[1]{\mathrel{\perp_{\!B,#1}^{\mathrm{mut}}}}
\newcommand{\bj}[1]{\mathrel{\perp_{\!B,#1}}}
\newcommand{\iso}[1]{\mathrel{\perp_{\!I,#1}}}
\newtheorem*{quotedtheorem}{Schlumprecht's theorem}
\newcommand{\Sn}[1]{S_{#1}}
\colorlet{blue}{black}
\AtBeginDocument{\hypersetup{hidelinks,allcolors=black}\color{black}}

\begin{document}
\title[The distortion problem]{The distortion problem and orthogonality}
\author[1]{\fnm{Zhiyao} \sur{Fang}}
\email{fzy73@stu.aqnu.edu.cn}
\author*[1]{\fnm{Qi} \sur{Liu}}
\email{liuq67@aqnu.edu.cn}
\affil[1]{\orgdiv{School of Mathematics and Statistics},
\orgname{Anqing Normal University},
\orgaddress{\city{Anqing},\postcode{246133},\state{Anhui},\country{China}}}
\abstract{Schlumprecht [Israel J. Math. 76 (1991), 81--95]
constructed an arbitrarily distortable Banach space.
We further study distortion on this space under orthogonality
constraints and obtain the following results.
\begin{enumerate}
\renewcommand{\labelenumi}{(\alph{enumi})}
\setlength{\itemsep}{0.3em}
\setlength{\parsep}{0pt}
\item The exact distortion of his renormings is preserved
in every infinite-dimensional subspace when the unit vectors
satisfy both mutual Birkhoff--James and isosceles orthogonality
for any finite family of equivalent norms.
\item There are basic sequences of unit vectors satisfying
both orthogonality conditions pairwise whose norms in the
renormed space converge to the lower distortion endpoint.
\end{enumerate}
The proof combines finite-dimensional norm estimates with
supporting functionals and the Borsuk--Ulam theorem. Further
results concern finite and approximate Roberts conditions
and an obstruction to simultaneous Pythagorean orthogonality
for auxiliary norms.}
\maketitle

\section{Introduction}\label{sec:introduction}

Let $N$ be a norm on a real vector space $X$. Birkhoff--James orthogonality is the relation
\begin{equation}\label{eq:bj-definition}
x\bj{N}y
\quad\Longleftrightarrow\quad
N(x+ty)\ge N(x)\quad\text{for every }t\in\R.
\end{equation}
It originates in the work of Birkhoff and James \cite{birkhoff,james}. For any norm $q$, we write $x\mut{q}y$ if both
$x\bj{q}y$ and $y\bj{q}x$ hold. For a subspace $F$, the notation $x\bj{N}F$ means that $N(x+u)\ge N(x)$ for every $u\in F$. Orthogonality to a subspace has this single simultaneous meaning throughout the paper.

Isosceles orthogonality is defined by
\begin{equation}\label{eq:isosceles-definition}
x\iso{N}y\quad\Longleftrightarrow\quad N(x+y)=N(x-y).
\end{equation}

For any norm $q$, the notation $x\bj{q}y$ means
$q(x+ty)\ge q(x)$ for every $t\in\R$, and $x\iso{q}y$
means $q(x+y)=q(x-y)$.
For a subspace $F$, we write $x\bj{q}F$ if
$q(x+u)\ge q(x)$ for every $u\in F$.

The two relations are compatible: in an inner product space,
every orthogonal pair satisfies both Birkhoff--James and
isosceles orthogonality.
Beyond the inner product setting, Lemmas~\ref{lem:isosceles-selection}
and~\ref{lem:low-selection} below provide simultaneous
orthogonality under suitable dimension and subspace conditions.

Isosceles orthogonality is symmetric. In the distortion results, the
witnesses are normalized with respect to $N$, while prescribed
auxiliary vectors retain their given lengths. All isosceles
equalities are imposed at these fixed scales.

For applications of Birkhoff--James and isosceles
orthogonality, as well as related notions of orthogonality,
see \cite{alonso-martini-wu,BhatiaSemrl1999,ZamaniMoslehian2016,
	MartinEtAl2024,ChengFelder2023,
	Alonso1993,Bottazzi2019,Cho1992,ChengRossSeco2022}.

For a second, equivalent norm $M$ on an infinite-dimensional Banach space $(X,N)$, define
\begin{equation}\label{eq:ordinary-distortion}
	D_M(Y)=\sup\left\{\frac{M(x)}{M(y)}:x,y\in\Sn{N}(Y)\right\},
	\qquad
	\Sn{N}(Y)=\{x\in Y:N(x)=1\}.
\end{equation}
The norm $M$ is a $\lambda$-distortion if $D_M(Y)\ge\lambda$
for every infinite-dimensional closed subspace $Y$ of $X$
\cite[p.~81]{schlumprecht}.
Schlumprecht \cite[Theorem~3]{schlumprecht} made a fundamental
contribution to distortion theory by constructing a Banach
space that admits such a renorming for every $\lambda>1$;
his construction subsequently became a key ingredient in
the solution of the classical distortion problem for $\ell_2$
by Odell and Schlumprecht \cite{odell-schlumprecht}.

To state the precise result used for comparison, let $S$ be the real Schlumprecht space with original norm $N$, put $f(t)=\log_2(t+1)$, and define
\[
\nu_\ell(x)=\frac1{f(\ell)}
\sup_{E_1<\cdots<E_\ell}\sum_{r=1}^{\ell}N(E_rx),\qquad \ell\ge2.
\]
Here $E_rx$ denotes coordinate restriction to a finite set,
and the sets are successive. The construction is recalled
in Section~\ref{sec:schlumprecht}. The following formulation
follows from Schlumprecht's Theorem~3 and the norm comparison
immediately preceding it \cite[p.~85]{schlumprecht}.

\begin{quotedtheorem}
	For each $\ell\ge2$, one has $N/f(\ell)\le\nu_\ell\le N$.
	For every infinite-dimensional subspace $Y\subset S$ and
	every $\varepsilon>0$, there exist $x,y\in\Sn{N}(Y)$ such that
	\begin{equation}\label{eq:schlumprecht-original}
		\nu_\ell(x)\ge1-\varepsilon,
		\qquad
		\nu_\ell(y)\le\frac{1+\varepsilon}{f(\ell)}.
	\end{equation}
\end{quotedtheorem}

The norm comparison gives $1/f(\ell)\le\nu_\ell(z)\le1$ for every
$z\in\Sn{N}(Y)$. Hence each ratio of two such norm values
is at most $f(\ell)$. For $0<\varepsilon<1$, the vectors
in \eqref{eq:schlumprecht-original} give a ratio at least
$f(\ell)(1-\varepsilon)/(1+\varepsilon)$. Letting
$\varepsilon\downarrow0$ gives the exact formula
\[
\sup\left\{\frac{\nu_\ell(x)}{\nu_\ell(y)}:
  x,y\in\Sn{N}(Y)\right\}=f(\ell).
\]

Subsequent work linked distortion to classical sequence spaces
and asymptotic geometry. Odell and Schlumprecht
\cite{odell-schlumprecht1993} announced their solution to
the distortion problem for Hilbert space and constructed
families of nearly orthogonal subsets of its unit sphere,
each meeting every infinite-dimensional closed subspace.
Their subsequent paper \cite{odell-schlumprecht} developed
the full argument establishing arbitrary distortability
of $\ell_2$, using Schlumprecht's space. Milman and Tomczak-Jaegermann \cite{milman-tomczak}
showed that bounded distortion forces the existence of a
subspace with asymptotic $\ell_p$ structure, including the
asymptotic $c_0$ case, through stabilization of norms.
Maurey \cite{maurey-distortion} proved that every Banach space
with an unconditional basis which does not contain
$\ell_1^n$ uniformly has an arbitrarily distortable subspace.
His argument transfers asymptotic sets from $\ell_2$ to
asymptotic $\ell_p$ spaces. Androulakis and Odell
\cite{androulakis-odell} established arbitrary distortability
of regular mixed Tsirelson spaces satisfying
$\theta_n/\theta^n\to0$, where
$\theta=\lim_{n\to\infty}\theta_n^{1/n}$, and calculated
associated asymptotic $\ell_1$ constants.

Orthogonality constraints restrict the pairs available to
witness distortion, so it is natural to ask whether they
reduce the hereditary distortion constant.
Proposition~\ref{prop:directed} shows that directed
Birkhoff--James orthogonality leaves this constant unchanged.
We then study whether mutual Birkhoff--James and isosceles
orthogonality can be imposed simultaneously.
On Schlumprecht's space, Corollary~\ref{cor:ell-exact}
preserves the exact distortion of his original renormings
under both conditions, showing that the full norm oscillation
can be realized by pairs satisfying these additional
geometric constraints.

Our main result concerns the real Schlumprecht space $S$, whose construction is recalled in Section~\ref{sec:schlumprecht}. The auxiliary norms may include $N$ and the distorting norm itself, and removing any orthogonality constraints preserves the same supremum; removing all of them recovers arbitrary distortability. For the original renormings $\nu_\ell$, Corollary~\ref{cor:ell-exact} retains the exact distortion $f(\ell)$ and the endpoint estimates in \eqref{eq:schlumprecht-original}. The proof extends the long-average argument of \cite[Theorem~3]{schlumprecht} to finite-dimensional spaces whose entire unit spheres lie near the lower renorming endpoint. The transfer principle in Section~\ref{sec:transfer} then imposes exact orthogonality after perturbation into the prescribed subspace, and interpolation with $N$ yields every $\lambda>1$. Further consequences include prescribed additional constraints and basic sequences approaching the lower endpoint; see ~\ref{cor:schlumprecht-sequence}.

Corollary~\ref{cor:ell-exact} preserves the supremum $f(\ell)$
for the same norm $\nu_\ell$ when the pairs in the formula
above are required to satisfy both orthogonality conditions
for each prescribed auxiliary norm. This gives
\eqref{eq:main-equality} at $\lambda=f(\ell)$ with
$M_\lambda=\nu_\ell$. For arbitrary $\lambda>1$, the
interpolation in \eqref{eq:interpolated-norm} gives the
required norm $M_\lambda$.
\begin{theorem}\label{thm:main}
	For every $\lambda>1$ there is an equivalent norm $M_\lambda$
	on $S$ such that
	\begin{equation}\label{eq:main-comparison}
		\lambda^{-1}N(x)\le M_\lambda(x)\le N(x)
		\qquad(x\in S),
	\end{equation}
	and, for every finite family $q_1,\ldots,q_m$ of equivalent
	norms on $S$ and every infinite-dimensional closed subspace
	$Y\subset S$,
	\begin{equation}\label{eq:main-equality}
		\sup\left\{
		\frac{M_\lambda(x)}{M_\lambda(y)}:
		x,y\in\Sn{N}(Y),\
		x\mut{q_j}y,\ x\iso{q_j}y
		\ (1\le j\le m)
		\right\}=\lambda.
	\end{equation}
	The family of norms can be specified after $M_\lambda$
	has been chosen.
\end{theorem}

\begingroup\color{black}
Section~\ref{sec:further-constraints} shows that excluding both
relations in the original norm also preserves the distortion
supremum, and examines additional Pythagorean and Roberts conditions.
\endgroup

All vector spaces and functionals in this paper are real.
Subspaces of a Banach space that are quantified as
infinite-dimensional are understood to be closed.
For a norm $q$, the dual norm of a functional is denoted by
$q^*$. For a linear operator $T$ on a normed space $(E,q)$,
we write
\[
\norm{T}_q
=\sup_{z\in E\setminus\{0\}}\frac{q(Tz)}{q(z)}.
\]
No smoothness is assumed for any of the norms.
\section{Simultaneous orthogonality for finitely many norms}\label{sec:selection}

We first record the support-functional form of orthogonality to a subspace. If $x\ne0$ and $x\bj{q}F$, then $x\notin F$. On $F+\R x$ the formula
\begin{equation}\label{eq:support-on-span}
p(u+tx)=tq(x)\qquad(u\in F,\ t\in\R)
\end{equation}
defines a functional of $q$-dual norm one. Indeed, for $t\ne0$,
\[
q(u+tx)=|t|q(x+u/t)\ge |t|q(x),
\]
and the inequality for $t=0$ is immediate. Hahn--Banach extends $p$ to the ambient space with $q^*(p)=1$, $p(x)=q(x)$ and $p|_F=0$. Conversely, these properties imply $x\bj{q}F$. This is the subspace form of James' characterization \cite{james}.

\begin{lemma}\label{lem:reverse-selection}
Let $q_1,\ldots,q_m$ be norms on a real vector space $V$, let $F\subset V$ have dimension $d$, and let $Z\subset V$ have dimension at least $md+1$. There exists $z\in Z\setminus\{0\}$ such that
\begin{equation}\label{eq:reverse-selection}
z\bj{q_j}F\qquad(1\le j\le m).
\end{equation}

\end{lemma}
\begin{proof}
If $d=0$, any nonzero vector of $Z$ works. Suppose $d\ge1$. Replace $Z$ by a subspace of dimension $md+1$, and restrict attention to the finite-dimensional space $E=F+Z$. Give $E$ an auxiliary Euclidean structure. Choose a nonnegative even function $\rho\in C_c^\infty(E)$ with integral one. For $\delta>0$ put
\[
h_{j,\delta}(v)=\int_E q_j(v-\delta w)\rho(w)\,dw.
\]
Each $h_{j,\delta}$ is smooth, convex and even. Also
\[
|h_{j,\delta}(v)-q_j(v)|
\le\delta\int_E q_j(w)\rho(w)\,dw,
\]
so the convergence to $q_j$ is uniform on $E$.

Choose a basis $u_1,\ldots,u_d$ of $F$. On the Euclidean unit sphere of $Z$, which is a sphere of dimension $md$, consider
\[
z\longmapsto
\big(Dh_{j,\delta}(z)[u_i]\big)_{
1\le j\le m,\,1\le i\le d}\in\R^{md}.
\]
The map is continuous and odd. The odd-map form of the Borsuk--Ulam theorem \cite[Chapter~2]{matousek} gives a zero $z_\delta$. Since the derivative of $h_{j,\delta}$ at $z_\delta$ vanishes on $F$, convexity gives
\begin{equation}\label{eq:smoothed-minimum}
h_{j,\delta}(z_\delta+u)\ge h_{j,\delta}(z_\delta)
\qquad(u\in F,\ 1\le j\le m).
\end{equation}
Along a sequence $\delta\downarrow0$, compactness gives $z_\delta\to z$ on the Euclidean unit sphere of $Z$. For an arbitrary fixed $u\in F$, uniform convergence and continuity allow passage to the limit in \eqref{eq:smoothed-minimum}. Thus $q_j(z+u)\ge q_j(z)$ for all $u$ and $j$. The vector $z$ is nonzero.

\end{proof}

\begin{lemma}\label{lem:isosceles-selection}
Let $N,q_1,\ldots,q_m$ be norms on a real vector space $V$, let $F\subset V$ have dimension $d$, and fix $u_1,\ldots,u_n\in V$. If $Z\subset V$ has dimension at least $m(d+n)+1$, there exists $z\in Z$ such that
\[
N(z)=1,\qquad z\bj{q_j}F,\qquad
q_j(z+u_i)=q_j(z-u_i)\quad(1\le i\le n,\ 1\le j\le m).
\]
The prescribed unit sphere is $\Sn{N}(Z)$; no condition $q_j(z)=1$ is imposed. If $u_i\in F$, the conclusion gives both $z\bj{q_j}u_i$ and $z\iso{q_j}u_i$. The opposite relations $u_i\bj{q_j}z$ require an additional kernel condition, as implemented in Lemma~\ref{lem:low-selection}. No optimality of this dimension bound is asserted.
\end{lemma}
\begin{proof}
Put $r=m(d+n)$, replace $Z$ by a subspace of dimension $r+1$, and work in the finite-dimensional space $E=F+Z+\spn\{u_1,\ldots,u_n\}$. On $E$ use the smooth even convex approximations $h_{j,\delta}$ constructed in Lemma~\ref{lem:reverse-selection}. Choose a basis $v_1,\ldots,v_d$ of $F$, omitting this list if $d=0$.

Let $\Sigma$ be the auxiliary Euclidean unit sphere of $Z$. The map
\[
R:\Sigma\longrightarrow\Sn{N}(Z),\qquad R(w)=w/N(w),
\]
is continuous, odd and onto. Thus normalization is built into the domain before the orthogonality equations are solved. For $z=R(w)$ consider the map whose coordinates are
\begin{equation}\label{eq:combined-odd-map}
\begin{gathered}
Dh_{j,\delta}(z)[v_h]\quad(1\le j\le m,\ 1\le h\le d),\\
h_{j,\delta}(z+u_i)-h_{j,\delta}(z-u_i)
\quad(1\le j\le m,\ 1\le i\le n).
\end{gathered}
\end{equation}
This is a continuous odd map from $\Sigma\cong S^r$ to $\R^r$. Oddness of the second group follows from evenness of $h_{j,\delta}$. The Borsuk--Ulam theorem gives a zero, with corresponding $N$-unit vector $z_\delta$. If $r=0$, the conditions are empty and any $N$-unit vector in $Z$ suffices.

Vanishing of the first group and convexity imply
\[
h_{j,\delta}(z_\delta+v)\ge h_{j,\delta}(z_\delta)\qquad(v\in F).
\]
Vanishing of the second group gives
$h_{j,\delta}(z_\delta+u_i)=h_{j,\delta}(z_\delta-u_i)$.
The $N$-unit sphere of $Z$ is compact, so along a sequence $\delta\downarrow0$ the vectors $z_\delta$ converge to an $N$-unit vector $z$. Uniform convergence of the approximations gives both sets of claimed relations. The argument takes place on the $N$-unit sphere from the outset; no separate rescaling is used after the isosceles equalities have been obtained.
\end{proof}

\begin{proposition}\label{prop:countable-obstruction}
Let $(X,N)$ be a nonzero separable real normed space and let $\eta>0$. There is a countable family of norms $(q_k)$ satisfying
\begin{equation}\label{eq:countable-comparison}
N\le q_k\le(1+\eta)N\qquad(k\in\N)
\end{equation}
such that no nonzero ordered pair $(x,y)$ satisfies $x\bj{q_k}y$ for every $k$.
\end{proposition}
\begin{proof}
Choose a dense sequence $(v_r)$ in $\Sn{N}(X)$ and functionals $f_r$ with $N^*(f_r)=1$ and $f_r(v_r)=1$. For $x\ne0$ and $v_r$ close to $x/N(x)$,
\[
f_r(x)/N(x)\ge1-N\big(x/N(x)-v_r\big).
\]
Thus $\sup_r f_r(x)=N(x)$ for every $x$.

Let $\mathcal H$ be the countable family of functionals
\[
h=(1-s)f_r+sf_t,
\qquad r,t\in\N,\quad s\in\mathbb Q\cap(0,1).
\]
Each has $N^*(h)\le1$. Put $c=1+\eta$ and define $q_h(x)=\max\{N(x),c|h(x)|\}$. This gives \eqref{eq:countable-comparison}.

Fix nonzero $x,y$. Choose $r$ with $f_r(x)>N(x)/c$, and choose $t$ with $f_t(y)\ne0$. For all sufficiently small positive $s$, the corresponding $h$ satisfies $h(x)>N(x)/c$. At most one such $s$ has $h(y)=0$, so a rational $s$ can be chosen with $h(y)\ne0$. Near $u=0$, strict dominance and continuity give
\[
q_h(x+uy)=c\big(h(x)+uh(y)\big).
\]
This function has nonzero derivative at zero and therefore does not attain a minimum there. Hence $x\not\bj{q_h}y$.
\end{proof}

\section{Distortion under orthogonality constraints}\label{sec:transfer}

We first show that imposing directed Birkhoff--James
orthogonality preserves hereditary distortion. The argument
combines the hereditary distortion hypothesis with
supporting-functional kernels and applies to every real
Banach space.

\begin{proposition}\label{prop:directed}
Let $M$ be a $\lambda$-distortion of a real Banach space $(X,N)$, and let $q_1,\ldots,q_m$ be equivalent norms on $X$. For every infinite-dimensional subspace $Y$,
\begin{equation}\label{eq:directed-distortion}
\sup\left\{\frac{M(x)}{M(y)}:
x,y\in\Sn{N}(Y),\ x\bj{q_j}y\ (1\le j\le m)\right\}\ge\lambda.
\end{equation}
\end{proposition}
\begin{proof}
For infinite-dimensional $Z\subset Y$, put $a(Z)=\inf_{\Sn{N}(Z)}M$ and $b(Z)=\sup_{\Sn{N}(Z)}M$, and let $a_* =\sup_{Z\subset Y}a(Z)$. Equivalence of $M$ and $N$ gives $0<a_*<\infty$, and distortion gives $b(Z)\ge\lambda a(Z)$.

Given sufficiently small $\delta>0$, choose $Z$ with $a(Z)>a_*-\delta$, and then $x\in\Sn{N}(Z)$ with
\[
M(x)>b(Z)-\delta>\lambda a_*-(\lambda+1)\delta.
\]
For each $j$ choose a $q_j$-norming functional at $x$ on $Z$. Their common kernel $W$ is infinite-dimensional and satisfies $a(W)\le a_*$. Choose $y\in\Sn{N}(W)$ with $M(y)<a_*+\delta$. James' characterization \cite[Corollary~2.2]{james} gives $x\bj{q_j}y$ for every $j$, while
\[
\frac{M(x)}{M(y)}>
\frac{\lambda a_*-(\lambda+1)\delta}{a_*+\delta}.
\]
Letting $\delta\downarrow0$ proves the assertion.
\end{proof}

Conversely, if \eqref{eq:directed-distortion} holds on every
infinite-dimensional subspace for a fixed finite family,
then $M$ is a $\lambda$-distortion, since removing the
orthogonality constraints can only increase the supremum.
Taking the family to be $\{N\}$ therefore gives an equivalent
formulation of hereditary $\lambda$-distortion using directed
orthogonality in the original norm. This equivalence concerns
the same lower bound on every infinite-dimensional subspace;
it does not assert equality of the constrained and
unconstrained suprema on each fixed subspace. The kernel
argument alone supplies no reverse orthogonality.

We now pass to mutual orthogonality. To keep the second vector near the lower endpoint while imposing the reverse conditions, we require uniformly small renormed lengths on a whole finite-dimensional sphere. Suppose $N,M$ are equivalent norms and $a>0$. The condition is
\begin{equation}\label{eq:lower-property}
\begin{gathered}
\text{for every infinite-dimensional }Y,\ d\in\N,\ \varepsilon>0,\\
\text{there is }E\subset Y,\quad \dim E=d,\quad
M(z)\le(a+\varepsilon)N(z)\quad(z\in E).
\end{gathered}
\end{equation}
Thus the entire unit sphere of $E$ lies near the lower endpoint.

\begin{lemma}\label{lem:low-selection}
Suppose \eqref{eq:lower-property} holds. Let $q_1,\ldots,q_m$ be equivalent norms on $X$, let $Y\subset X$ be infinite-dimensional, and let $u_1,\ldots,u_n\in X$. For every $\varepsilon>0$ there is $y\in\Sn{N}(Y)$ such that
\begin{equation}\label{eq:low-selection-conclusion}
\begin{gathered}
M(y)\le a+\varepsilon,\\
y\iso{q_j}u_i,\quad u_i\bj{q_j}y
\quad(1\le i\le n,\ 1\le j\le m),\\
y\bj{q_j}\spn\{u_1,\ldots,u_n\}\quad(1\le j\le m).
\end{gathered}
\end{equation}
For each $j$, a $q_j$-norming functional at $y$ may be chosen to annihilate the whole indicated span.
\end{lemma}
\begin{proof}
Discard zero members of the finite list, which impose automatic conditions, and relabel the remaining $n$ vectors. Write $F=\spn\{u_1,\ldots,u_n\}$ and $d=\dim F$. For each $i,j$, choose a functional $p_{j,i}$ with $q_j^*(p_{j,i})=1$ and $p_{j,i}(u_i)=q_j(u_i)$, and put
\[
W=Y\cap\bigcap_{j=1}^m\bigcap_{i=1}^n\ker p_{j,i}.
\]
This subspace is infinite-dimensional. Every $z\in W$ satisfies $u_i\bj{q_j}z$, since
\[
q_j(u_i+tz)\ge p_{j,i}(u_i+tz)=q_j(u_i)
\qquad(t\in\R).
\] By \eqref{eq:lower-property}, choose $Z\subset W$ with $\dim Z=m(d+n)+1$ and $M(z)\le(a+\varepsilon)N(z)$ for all $z\in Z$.

Apply Lemma~\ref{lem:isosceles-selection} to obtain $y\in\Sn{N}(Z)$ which is $q_j$-orthogonal to $F$ and isosceles orthogonal to each $u_i$ for all $j$. Because $u_i\in F$, the lemma gives $y\bj{q_j}u_i$, while $y\in W$ supplies $u_i\bj{q_j}y$. Hence the same $N$-unit vector satisfies mutual Birkhoff--James and isosceles orthogonality for every $i,j$. Its normalization is already fixed, and the inequality on $Z$ gives the estimate for $M(y)$. Formula~\eqref{eq:support-on-span} and Hahn--Banach give the final assertion.
\end{proof}

\begin{theorem}\label{thm:endpoint-transfer}
Let $(X,N)$ be a real Banach space, let $0<a\le b<\infty$, and let $M$ be a norm such that
\begin{equation}\label{eq:global-endpoints}
aN(x)\le M(x)\le bN(x)\qquad(x\in X).
\end{equation}
Assume \eqref{eq:lower-property} and
\begin{equation}\label{eq:upper-property}
\sup_{x\in\Sn{N}(Y)}M(x)=b
\quad\text{for every infinite-dimensional }Y\subset X.
\end{equation}
For any finite family of equivalent norms $q_1,\ldots,q_m$, any infinite-dimensional $Y$, and any finite list $u_1,\ldots,u_n\in X$, the supremum of $M(x)/M(y)$ over $x,y\in\Sn{N}(Y)$ subject to
\begin{equation}\label{eq:transfer-constraints}
\begin{gathered}
x\mut{q_j}y,\quad x\iso{q_j}y,\qquad
y\bj{q_j}\spn\{x,u_1,\ldots,u_n\}\quad(1\le j\le m),\\
u_i\bj{q_j}y,\quad u_i\iso{q_j}y
\quad(1\le i\le n,\ 1\le j\le m)
\end{gathered}
\end{equation}
is exactly $b/a$.
\end{theorem}
\begin{proof}
The upper bound follows from \eqref{eq:global-endpoints}. For $0<\varepsilon<b$, choose $x\in\Sn{N}(Y)$ with $M(x)>b-\varepsilon$, using \eqref{eq:upper-property}. Apply Lemma~\ref{lem:low-selection} to the list $x,u_1,\ldots,u_n$ to obtain $y\in\Sn{N}(Y)$ satisfying all the constraints and $M(y)\le a+\varepsilon$. Then
\[
\frac{M(x)}{M(y)}>\frac{b-\varepsilon}{a+\varepsilon}.
\]
Let $\varepsilon\downarrow0$.
\end{proof}

The same selection principle gives infinite basic sequences with exact coordinate-projection estimates.

\begin{theorem}\label{thm:basic-sequence}
Assume \eqref{eq:lower-property} and $M\ge aN$. Let $Y$ be infinite-dimensional, $x_0\in\Sn{N}(Y)$, $\eta>0$, and let $(\varepsilon_k)_{k\ge1}$ be positive. For any finite family of equivalent norms, there is a sequence $(x_k)_{k\ge0}\subset\Sn{N}(Y)$ with the prescribed first term, such that
\begin{equation}\label{eq:sequence-low}
a\le M(x_k)\le a+\varepsilon_k\qquad(k\ge1).
\end{equation}
It is a basic sequence for each prescribed norm and for $N$, with basis constant at most $1+\eta$. Every coordinate projection has norm one for all these norms. The vectors at this $N$-unit normalization are pairwise both mutually Birkhoff--James orthogonal and isosceles orthogonal for each of the norms.
\end{theorem}
\begin{proof}
Adjoin $N$ to the finite family if necessary, and denote the resulting norms by $q_1,\ldots,q_m$. We construct $x_k$, norming functionals $p_{j,k}$, and finite sets of functionals $\mathcal A_{j,k}$.

After $x_0,\ldots,x_k$ have been chosen, write $F_k=\spn\{x_0,\ldots,x_k\}$. Compactness of the $q_j$-unit sphere of $F_k$ and Hahn--Banach give a finite set $\mathcal A_{j,k}$ of $q_j$-dual unit functionals on $X$ such that
\begin{equation}\label{eq:finite-norming}
\max_{p\in\mathcal A_{j,k}}|p(v)|\ge\frac{q_j(v)}{1+\eta}
\qquad(v\in F_k).
\end{equation}
To see this, take a norming functional at each point of that sphere. The open sets where its absolute value exceeds $(1+\eta)^{-1}$ cover the sphere; a finite subcover suffices.

At the initial stage choose $p_{j,0}$ norming $x_0$. Suppose the data have been chosen through stage $k$. Intersect $Y$ with the kernels of every $p_{j,h}$ and every functional in $\mathcal A_{j,h}$, for $0\le h\le k$ and $1\le j\le m$. This is an infinite-dimensional subspace $W_k$. Lemma~\ref{lem:low-selection}, applied within $W_k$ to the list $x_0,\ldots,x_k$, gives an $N$-unit vector $x_{k+1}$ with
\[
M(x_{k+1})\le a+\varepsilon_{k+1},\quad
x_{k+1}\iso{q_j}x_h\ (0\le h\le k),\qquad
x_{k+1}\bj{q_j}F_k\quad(1\le j\le m).
\]
Choose $p_{j,k+1}$ with $q_j^*(p_{j,k+1})=1$, $p_{j,k+1}(x_{k+1})=q_j(x_{k+1})$ and $p_{j,k+1}|_{F_k}=0$, using \eqref{eq:support-on-span}. Then choose the sets $\mathcal A_{j,k+1}$. Induction gives
\begin{equation}\label{eq:biorthogonal-supports}
p_{j,k}(x_h)=0\quad(h\ne k),\qquad
p_{j,k}(x_k)=q_j(x_k).
\end{equation}
In particular the vectors are linearly independent. Their pairwise isosceles equalities were imposed at each selection stage and remain unchanged thereafter.

Let $v\in F_k$ and let $w$ be a finite combination of $x_h$ with $h>k$. All functionals in $\mathcal A_{j,k}$ annihilate $w$. Therefore \eqref{eq:finite-norming} gives
\[
q_j(v+w)\ge q_j(v)/(1+\eta).
\]
The initial-segment projections $P_k$ on the algebraic span consequently satisfy $\norm{P_k}_{q_j}\le1+\eta$. They extend to its closure $E$. For $z\in E$, approximate $z$ by a finite combination $v$; when $k$ contains its support,
\[
q_j(P_kz-z)\le(2+\eta)q_j(z-v).
\]
Hence $P_kz\to z$, proving the basic-sequence assertion.

Let $Q_k$ be the $k$th coordinate projection. On the algebraic span, \eqref{eq:biorthogonal-supports} yields
\[
q_j(Q_kz)=|p_{j,k}(z)|\le q_j(z).
\]
This extends to $E$; equality holds on $x_k$, so $\norm{Q_k}_{q_j}=1$. Applying this inequality to $x_k+tx_h$ proves $x_k\bj{q_j}x_h$ for $h\ne k$. Interchanging $k,h$ gives mutual orthogonality. The lower inequality in \eqref{eq:sequence-low} follows from $M\ge aN$.
\end{proof}

For each $q_j$, normalization of the vectors in Theorem~\ref{thm:basic-sequence} gives an Auerbach basis of their closed span: the normalized vectors and their coefficient functionals both have norm one. The common objects across the norms are the coordinate directions and projections.

\section{The Schlumprecht space}\label{sec:schlumprecht}

Set $f(t)=\log_2(t+1)$ for $t\ge1$. The real Schlumprecht space $S$ is the completion of $c_{00}$ under the implicit norm
\begin{equation}\label{eq:schlumprecht-norm}
	N(x)=\max\left\{
	\norm{x}_\infty,\,
	\sup_{\substack{r\ge2\\E_1<\cdots<E_r}}
	\frac{1}{f(r)}\sum_{i=1}^{r}N(E_i x)
	\right\}.
\end{equation}
Here $E_i$ are finite nonempty subsets of $\N$, $E<F$ means $\max E<\min F$, and $Ex$ is the coordinate restriction of $x$ to $E$. The canonical basis is normalized and $1$-unconditional \cite[Proposition~2]{schlumprecht}. In particular all coordinate projections are contractive for $N$.

For $\ell\ge2$ define
\begin{equation}\label{eq:ell-norm}
\nu_\ell(x)=\sup_{E_1<\cdots<E_\ell}
\frac1{f(\ell)}\sum_{r=1}^{\ell}N(E_rx).
\end{equation}
These norms satisfy
\begin{equation}\label{eq:ell-comparison}
f(\ell)^{-1}N(x)\le\nu_\ell(x)\le N(x)\qquad(x\in S).
\end{equation}
We use two results from \cite{schlumprecht}. Its Lemma~5 says that $\ell_1$ is finitely block represented in every infinite block subspace: for each $L\in\N$ and $\delta>0$ there are normalized successive block vectors $v_1,\ldots,v_L$ satisfying
\begin{equation}\label{eq:l1-blocks}
\frac1{1+\delta}\sum_{i=1}^L|c_i|
\le N\left(\sum_{i=1}^L c_iv_i\right)
\le\sum_{i=1}^L|c_i|\qquad(c_i\in\R).
\end{equation}
Its Theorem~3 gives, in every infinite-dimensional $Y\subset S$ and for every $\varepsilon>0$, an $x\in\Sn{N}(Y)$ with $\nu_\ell(x)\ge1-\varepsilon$. Thus
\begin{equation}\label{eq:ell-upper-endpoint}
\sup_{x\in\Sn{N}(Y)}\nu_\ell(x)=1.
\end{equation}
The next estimate supplies the required uniform finite-dimensional lower endpoint.

\begin{lemma}\label{lem:average-space}
Fix $\ell,d,k\in\N$ with $\ell\ge2$ and $\delta>0$. Suppose $v_1,\ldots,v_{dk}$ are normalized successive block vectors satisfying \eqref{eq:l1-blocks}. Put
\[
z_h=\frac1k\sum_{i=(h-1)k+1}^{hk}v_i\qquad(1\le h\le d).
\]
Then for every $z\in\spn\{z_1,\ldots,z_d\}$,
\begin{equation}\label{eq:average-bound}
\nu_\ell(z)\le
\frac{(1+\delta)(1+2\ell/k)}{f(\ell)}N(z).
\end{equation}
\end{lemma}
\begin{proof}
Write $z=\sum_{h=1}^d a_hz_h$ and $A=\sum_h|a_h|$. If $A=0$ there is nothing to prove. By \eqref{eq:l1-blocks},
\begin{equation}\label{eq:average-lower}
N(z)\ge A/(1+\delta).
\end{equation}
Fix $E_1<\cdots<E_\ell$ in \eqref{eq:ell-norm}. Enlarge each $E_r$ to its interval hull. Their ordering is preserved and each $N(E_rz)$ can only increase, by contractivity of coordinate projections. We may therefore estimate successive intervals.

For any interval $E_r$, at most two of the successive blocks $v_i$ have support that meets $E_r$ without being contained in it: any such block must cross one of its endpoints. For blocks wholly contained in $E_r$, use the triangle inequality. Summing their contributions over $r$ gives at most $A$, since no nonempty block support can be contained in two of the disjoint intervals. Every remaining contribution has norm at most $A/k$: its coefficient has absolute value $|a_h|/k\le A/k$, and $N(E_rv_i)\le1$. There are at most $2\ell$ such contributions. Hence
\[
\sum_{r=1}^{\ell}N(E_rz)\le(1+2\ell/k)A.
\]
Taking the supremum in \eqref{eq:ell-norm} and using \eqref{eq:average-lower} proves \eqref{eq:average-bound}.
\end{proof}

\begin{proposition}\label{prop:lower-endpoint}
For every infinite-dimensional $Y\subset S$, every $d\in\N$, every $\ell\ge2$ and every $\varepsilon>0$, there is a $d$-dimensional subspace $E\subset Y$ such that
\begin{equation}\label{eq:uniform-lower-endpoint}
\nu_\ell(z)\le\frac{1+\varepsilon}{f(\ell)}N(z)
\qquad(z\in E).
\end{equation}
\end{proposition}
\begin{proof}
If $Y$ is an infinite block subspace, choose $k$ large and then $\delta>0$ small enough that $(1+\delta)(1+2\ell/k)\le1+\varepsilon$. Choose $dk$ blocks by \eqref{eq:l1-blocks} and apply Lemma~\ref{lem:average-space}. Their $d$ averages are nonzero successive blocks and thus linearly independent.

We give the perturbation argument for a general $Y$. Put $a=f(\ell)^{-1}$ and choose
\begin{equation}\label{eq:perturbation-size}
0<\tau<\frac{\varepsilon a}{2\{1+(1+\varepsilon)a\}}.
\end{equation}
Choose positive numbers $\eta_i$ with $\sum_i\eta_i<\tau$. We recursively construct $N$-unit vectors $y_i\in Y$, normalized successive block vectors $b_i$, and integers $0=p_0<p_1<\cdots$ such that $N(y_i-b_i)<\eta_i$.

Having chosen $p_{i-1}$, select $y_i\in\Sn{N}(Y)$ with its first $p_{i-1}$ coordinates zero. This is possible because the corresponding coordinate kernel has finite codimension in $Y$. Choose $p_i>p_{i-1}$ so that
\[
N\big(y_i-P_{(p_{i-1},p_i]}y_i\big)<\eta_i/2,
\]
where $P_I$ denotes coordinate restriction to $I$. The truncation is nonzero. Normalize it to obtain $b_i$. The reverse triangle inequality gives $N(y_i-b_i)<\eta_i$.

For a finite combination $b=\sum_i c_ib_i$, contractivity of the coordinate projection onto the support of $b_i$ gives $|c_i|\le N(b)$. The map $Tb_i=y_i$ consequently satisfies
\begin{equation}\label{eq:perturbation-map}
N(Tb-b)\le\sum_i|c_i|N(y_i-b_i)\le\tau N(b).
\end{equation}
It extends to the closed block span $B=[b_i]$, is bounded below by $1-\tau$, and has range in $Y$.

Apply the block case to $B$, with error $\varepsilon/2$, to find a $d$-dimensional $E_0\subset B$ on which $\nu_\ell(b)\le a(1+\varepsilon/2)N(b)$. By \eqref{eq:ell-comparison} and \eqref{eq:perturbation-map}, for nonzero $b\in E_0$,
\[
\frac{\nu_\ell(Tb)}{N(Tb)}
\le\frac{a(1+\varepsilon/2)+\tau}{1-\tau}
\le a(1+\varepsilon).
\]
The last inequality is exactly the consequence of \eqref{eq:perturbation-size} obtained by multiplying by $1-\tau$. The space $E=T(E_0)$ has dimension $d$ and lies in $Y$.
\end{proof}

Schlumprecht's Theorem~3, together with the norm comparison
preceding it \cite[p.~85]{schlumprecht}, gives
$D_{\nu_\ell}(Y)=f(\ell)$ for every infinite-dimensional
subspace $Y\subset S$. The following corollary preserves
this exact value for the same norm $\nu_\ell$ while imposing
simultaneous mutual Birkhoff--James and isosceles orthogonality
for any finite family of equivalent norms.

\begin{corollary}\label{cor:ell-exact}
For every $\ell\ge2$, every finite family of equivalent norms $q_1,\ldots,q_m$ on $S$, and every infinite-dimensional $Y\subset S$,
\begin{equation}\label{eq:ell-exact}
\sup\left\{
\frac{\nu_\ell(x)}{\nu_\ell(y)}:
x,y\in\Sn{N}(Y),\
x\mut{q_j}y,\ x\iso{q_j}y\ (1\le j\le m)
\right\}=f(\ell).
\end{equation}
\end{corollary}
\begin{proof}
Apply Theorem~\ref{thm:endpoint-transfer} with $M=\nu_\ell$, $a=f(\ell)^{-1}$ and $b=1$. The hypotheses follow from \eqref{eq:ell-comparison}, \eqref{eq:ell-upper-endpoint} and Proposition~\ref{prop:lower-endpoint}.
\end{proof}

\begin{proof}[Proof of Theorem~\ref{thm:main}]
Choose $\ell\ge2$ with $f(\ell)\ge\lambda$, and put
\begin{equation}\label{eq:interpolated-norm}
\alpha=\frac{1-\lambda^{-1}}{1-f(\ell)^{-1}},
\qquad M_\lambda=(1-\alpha)N+\alpha\nu_\ell.
\end{equation}
Then $0<\alpha\le1$, so $M_\lambda$ is a norm. Formula~\eqref{eq:ell-comparison} gives
\[
\big(1-\alpha+\alpha/f(\ell)\big)N
=\lambda^{-1}N\le M_\lambda\le N.
\]
Equation~\eqref{eq:ell-upper-endpoint} gives the upper endpoint property with $b=1$. Proposition~\ref{prop:lower-endpoint} gives \eqref{eq:lower-property} with $a=\lambda^{-1}$: to obtain an additive error $\varepsilon$, apply that proposition with multiplicative error $\varepsilon f(\ell)/\alpha$. Theorem~\ref{thm:endpoint-transfer} now proves \eqref{eq:main-equality}.
\end{proof}

Theorem~\ref{thm:endpoint-transfer} also allows the lower
distortion witness to satisfy finitely many prescribed
additional orthogonality conditions.

\begin{corollary}\label{cor:schlumprecht-sequence}
Fix $\lambda>1$, a finite family of equivalent norms on $S$, an infinite-dimensional $Y\subset S$, and $\eta,\varepsilon>0$. There is an $N$-normalized basic sequence $(x_k)_{k\ge0}$ in $Y$ with
\[
M_\lambda(x_0)>1-\varepsilon,
\qquad\lim_{k\to\infty}M_\lambda(x_k)=\lambda^{-1},
\]
whose basis constant is at most $1+\eta$ and whose coordinate projections have norm one for every prescribed norm, for $N$, and for $M_\lambda$. At their $N$-unit normalization the vectors are pairwise isosceles orthogonal for all these norms as well.
\end{corollary}
\begin{proof}
Choose $x_0\in\Sn{N}(Y)$ with $M_\lambda(x_0)>1-\varepsilon$. Adjoin $N$ and $M_\lambda$ to the prescribed family and apply Theorem~\ref{thm:basic-sequence} with a positive sequence $\varepsilon_k\downarrow0$ and $a=\lambda^{-1}$.
\end{proof}

Since $S$ is separable, Proposition~\ref{prop:countable-obstruction} applies to its original norm. Thus the finite family in Theorem~\ref{thm:main} cannot be replaced by an arbitrary countable family, even if every norm in that family is within any prescribed multiplicative factor $1+\eta$ of $N$.

\begingroup
\color{black}
\section{Further orthogonality constraints}\label{sec:further-constraints}

Theorem~\ref{thm:main} preserves distortion under simultaneous
mutual Birkhoff--James and isosceles orthogonality. We now examine
what happens when both relations are excluded in the original norm,
or when an additional orthogonality relation is imposed.

\medskip\noindent\textbf{Excluding both orthogonality relations.}
We use the classical fact that on the unit circle of a real normed plane, each unit vector has exactly two isosceles-orthogonal unit vectors, which are negatives of one another. This is recorded, with attribution to Alonso~\cite{Alonso1994}.

\begin{lemma}\label{lem:density}
Let $(Y,N)$ be a real normed space with $\dim Y\ge2$. The set
\[
\mathcal G_N(Y)=\left\{(x,y)\in\Sn{N}(Y)^2:
\begin{array}{l}
x\not\bj{N}y,\quad y\not\bj{N}x,\\
N(x+y)\ne N(x-y)
\end{array}\right\}
\]
is open and dense in $\Sn{N}(Y)^2$.
\end{lemma}
\begin{proof}
First suppose that $Y=E$ is two-dimensional, and write $C=\Sn{N}(E)$. The subsets
\[
\begin{aligned}
\mathcal B&=\{(x,y)\in C^2:x\bj{N}y\},\\
\mathcal B^{\mathrm{op}}&=\{(x,y):(y,x)\in\mathcal B\},\\
\mathcal J&=\{(x,y)\in C^2:N(x+y)=N(x-y)\}
\end{aligned}
\]
are closed. For $\mathcal B$, this follows by intersecting, over all $t\in\R$, the closed conditions $N(x+ty)\ge1$; the other two cases follow from coordinate interchange and continuity.

Smooth points are dense in $C$. Indeed, at a nonsmooth point the cone of supporting functionals contains a nonempty open interval of Euclidean normal directions. These intervals have disjoint interiors for distinct points: a normal in the interior of such a cone exposes that point alone. Assigning a rational angle to each interval shows that there are at most countably many nonsmooth points. Every open arc of $C$ is uncountable, proving the asserted density.

Let $U,V$ be nonempty relatively open subsets of $C$. Choose a smooth point $x\in U$, and let $f_x$ be its unique norming functional. Then
\[
x\bj{N}y\quad\Longleftrightarrow\quad f_x(y)=0.
\]
For the reverse implication use $N(x+ty)\ge f_x(x+ty)=1$. For the forward implication, Hahn--Banach gives a norming functional at $x$ vanishing at $y$, and uniqueness identifies it with $f_x$. The line $\ker f_x$ meets $C$ in exactly two points. Hence $V$ contains $y$ with $(x,y)\notin\mathcal B$. Thus $\mathcal B$ has empty interior; interchanging the coordinates proves the same for $\mathcal B^{\mathrm{op}}$.

For each fixed $x\in U$, the planar isosceles-uniqueness result says that the section $\{y\in C:(x,y)\in\mathcal J\}$ consists of exactly two points. Therefore $\mathcal J$ also has empty interior. The intersection of the three open dense complements is open and dense in $C^2$.

For arbitrary $Y$, take any pair in $\Sn{N}(Y)^2$ and a two-dimensional subspace containing that pair, enlarging its span if necessary. Approximation in this plane gives approximation in $Y$. The orthogonality conditions on two given vectors depend only on their span, so they have the same meaning in the plane and in $Y$. Openness follows from the closedness arguments above, valid in every dimension.
\end{proof}

\begin{theorem}\label{thm:nonorthogonal}
Let $(X,N)$ be a real normed space. For every equivalent norm $M$ on $X$ and every subspace $Y\subset X$ with $\dim Y\ge2$,
\[
\sup_{(x,y)\in\mathcal G_N(Y)}\frac{M(x)}{M(y)}=D_M(Y).
\]
\end{theorem}
\begin{proof}
The ratio is continuous on $\Sn{N}(Y)^2$, since $M(y)>0$. Every pair can be approximated by pairs in $\mathcal G_N(Y)$ by Lemma~\ref{lem:density}, and their ratios converge to its ratio. This proves the lower bound; inclusion gives the upper bound.
\end{proof}

\begin{corollary}\label{cor:schlumprecht-nonorthogonal}
For the norm $M_\lambda$ supplied by Theorem~\ref{thm:main}, every infinite-dimensional closed $Y\subset S$ satisfies
\[
\sup_{\substack{x,y\in\Sn{N}(Y)\\
x\not\bj{N}y,\ y\not\bj{N}x,\ x\not\iso{N}y}}
\frac{M_\lambda(x)}{M_\lambda(y)}=\lambda.
\]
\end{corollary}
\begin{proof}
Theorem~\ref{thm:main} and the bounds $N/\lambda\le M_\lambda\le N$ give $D_{M_\lambda}(Y)=\lambda$. Apply Theorem~\ref{thm:nonorthogonal}.
\end{proof}
This is stronger than merely excluding mutual Birkhoff--James orthogonality: both directed relations are excluded separately. The assertion concerns orthogonality in $N$. It does not assert the analogous density when isosceles orthogonality is measured in an arbitrary auxiliary norm while normalization remains in $N$.

\medskip\noindent\textbf{A Pythagorean obstruction.}
We use James's definition of Pythagorean orthogonality
\cite{James1945}; see also \cite[p.~450]{Partington}:
\[
x\pyo{q}y\quad\Longleftrightarrow\quad q(x-y)^2=q(x)^2+q(y)^2.
\]
If isosceles orthogonality also holds, the same equality holds with $x+y$.

Among generalized orthogonality relations, Pythagorean
orthogonality presents distinctive difficulties, as
illustrated by questions in the theory of functional
equations. Szab\'o \cite{Szabo1997} proved that odd Pythagorean
orthogonally additive mappings into an Abelian group are additive
when the domain has dimension at least three. Sikorska
\cite[Problem~2.2]{Sikorska2015} posed the problems of
determining even Pythagorean orthogonally additive
mappings and odd Pythagorean orthogonally additive
mappings in dimension two.

\begin{theorem}\label{thm:pythagorean}
Let $N,M$ be norms and let $x,y\ne0$. If $x\pyo{q}y$ for all three norms $q=N,M,N+M$, then
\[
N(x)M(y)=N(y)M(x).
\]
In particular, if $N(x)=N(y)=1$, then $M(x)=M(y)$.
\end{theorem}
\begin{proof}
Put $a=N(x)$, $b=N(y)$, $c=M(x)$ and $d=M(y)$, all positive. The identities for $N$ and $M$ give
\[
N(x-y)=\sqrt{a^2+b^2},\qquad M(x-y)=\sqrt{c^2+d^2}.
\]
The identity for $N+M$ is consequently
\[
\bigl(\sqrt{a^2+b^2}+\sqrt{c^2+d^2}\bigr)^2
=(a+c)^2+(b+d)^2.
\]
After cancellation this becomes
\[
\sqrt{(a^2+b^2)(c^2+d^2)}=ac+bd.
\]
Both sides are nonnegative. Squaring and subtracting yields $(ad-bc)^2=0$, so $ad=bc$.
\end{proof}

\begin{corollary}\label{cor:failure}
Theorem~\ref{thm:main} cannot be strengthened by additionally requiring Pythagorean orthogonality for every prescribed auxiliary norm.
\end{corollary}
\begin{proof}
For any proposed $M_\lambda$, choose the permitted auxiliary family $(N,M_\lambda,N+M_\lambda)$. All three norms are equivalent. Every admissible $N$-unit pair would have ratio $M_\lambda(x)/M_\lambda(y)=1$ by Theorem~\ref{thm:pythagorean}. Thus the ratio supremum cannot equal $\lambda>1$: it is one if the admissible set is nonempty, and otherwise there are no witnesses.
\end{proof}
This obstruction concerns the simultaneous requirement for arbitrary finite auxiliary families. It does not rule out Pythagorean orthogonality for $N$ alone. The proof does not require the Birkhoff--James or isosceles conditions.

\medskip\noindent\textbf{Finite and approximate Roberts conditions.}
Roberts orthogonality \cite{Roberts1934} for a norm $q$
is defined by
\[
q(x+ty)=q(x-ty)\qquad(t\in\R).
\]
It implies isosceles orthogonality by setting $t=1$, and $x\bj{q}y$ because
\[
q(x)\le\tfrac12\bigl(q(x+ty)+q(x-ty)\bigr)=q(x+ty).
\]
Homogeneity and evenness turn the same equalities into $q(y+sx)=q(y-sx)$ for every $s$, giving $y\bj{q}x$ as well.

For geometric characterizations of Roberts orthogonality in
normed planes and a quantitative comparison with Birkhoff
orthogonality, see \cite[Sects.~3--5]{BalestroMartiniTeixeira2016}.

For the rest of this section, $M_\lambda$ is the norm defined in
\eqref{eq:interpolated-norm}. The proof of Theorem~\ref{thm:main}
establishes \eqref{eq:lower-property} and \eqref{eq:upper-property}
with $a=\lambda^{-1}$ and $b=1$. These endpoint properties and
Lemma~\ref{lem:low-selection} allow additional equalities at
finitely many prescribed parameters.

\begin{theorem}\label{thm:finite-roberts}
Fix $\lambda>1$, a finite family $q_1,\ldots,q_m$ of equivalent norms on $S$, and a finite set $T\subset\R$. In every infinite-dimensional closed $Y\subset S$, the supremum of $M_\lambda(x)/M_\lambda(y)$ over $x,y\in\Sn{N}(Y)$ remains $\lambda$ when, besides $x\mut{q_j}y$ and $x\iso{q_j}y$ for $1\le j\le m$, one requires
\[
q_j(x+ty)=q_j(x-ty)\quad(t\in T,\ 1\le j\le m).
\]
The finite set $T$ and the auxiliary norms can be specified after $M_\lambda$ has been chosen.
\end{theorem}
\begin{proof}
Put $a=\lambda^{-1}$ and choose $x\in\Sn{N}(Y)$ with $M_\lambda(x)>1-\epsilon$. Apply Lemma~\ref{lem:low-selection} to the finite list consisting of $x$ and all vectors $x/t$ with $t\in T\setminus\{0\}$. It gives $y\in\Sn{N}(Y)$ with $M_\lambda(y)\le a+\epsilon$ and the required mutual conditions for $x$.
For $t\ne0$, the imposed equality
\[
q_j(y+x/t)=q_j(y-x/t)
\]
multiplied by $|t|$, together with evenness, gives the desired equality at $t$. The equality at zero is automatic. The ratio exceeds $(1-\epsilon)/(a+\epsilon)$. Let $\epsilon\downarrow0$ and use the global upper bound $\lambda$.
\end{proof}

Zamani and Moslehian \cite[Sect.~2]{ZamaniMoslehian2015Roberts}
studied two notions of approximate Roberts orthogonality. One
uses $q(x+ty)+q(x-ty)$ to scale the difference of the two norm
values. In the next corollary we use $q(x)+|t|q(y)$ and ask
whether this uniform bound preserves the distortion supremum.

\begin{corollary}\label{cor:approximate-roberts}
Fix $\lambda>1$, $\eta>0$, and a finite family $q_1,\ldots,q_m$ of equivalent norms on $S$. Adding the uniform inequalities
\begin{equation}\label{eq:approximate-roberts}
\begin{gathered}
|q_j(x+ty)-q_j(x-ty)|
\le\eta\bigl(q_j(x)+|t|q_j(y)\bigr)\\
(t\in\R,\ 1\le j\le m).
\end{gathered}
\end{equation}
to the exact conditions $x\mut{q_j}y$ and $x\iso{q_j}y$ for $1\le j\le m$ leaves the supremum of $M_\lambda(x)/M_\lambda(y)$ over $x,y\in\Sn{N}(Y)$ equal to $\lambda$ in every infinite-dimensional closed $Y\subset S$.
\end{corollary}
\begin{proof}
Choose $0<c\le C$ with $cN\le q_j\le CN$ for every $j$. Choose $R>2C/(\eta c)$ and $h<\eta c/(2C)$. Take a finite $h$-net $T$ of $[-R,R]$ and impose all equalities from Theorem~\ref{thm:finite-roberts}.
For $g_j(t)=q_j(x+ty)-q_j(x-ty)$, the triangle inequality gives
\[
|g_j(t)-g_j(s)|\le2q_j(y)|t-s|\le2C|t-s|.
\]
For $|t|\le R$, choose $s\in T$ with $|t-s|\le h$. Since $g_j(s)=0$, we have $|g_j(t)|\le2Ch<\eta c\le\eta(q_j(x)+|t|q_j(y))$.
For $|t|>R$, evenness and the reverse triangle inequality give
\[
\begin{aligned}
|g_j(t)|&=|q_j(ty+x)-q_j(ty-x)|\le2q_j(x)\le2C\\
&<\eta c|t|\le\eta(q_j(x)+|t|q_j(y)).
\end{aligned}
\]
Every pair supplied with these finite equalities satisfies \eqref{eq:approximate-roberts}. Their ratio supremum is $\lambda$ by Theorem~\ref{thm:finite-roberts}, and no admissible ratio exceeds $\lambda$.
\end{proof}

Whether the full exact Roberts relations can be imposed for every
finite auxiliary family while retaining the supremum $\lambda$
remains open here. The pairs supplied for finite parameter sets
may vary with those sets, and the preceding argument provides no
compactness principle producing a common pair for all parameters.
\endgroup

The preceding results suggest introducing geometric constants
associated with the simultaneous orthogonality constraints.

\medskip
\noindent\textbf{Question.}
Which geometric properties of Banach spaces can be characterized
by simultaneous Birkhoff--James and isosceles orthogonality,
and how does its characterizing power compare with that of
either relation alone?
For related results, see
\cite{alonso-martini-wu,LiuWangBao}.

For an equivalent norm $|\cdot|$ on $(X,\|\cdot\|)$,
the distortion quantity (see \cite[p.~127]{androulakis-odell}) is
\[
D(Y,|\cdot|)
=
\sup\left\{
\frac{|y|}{|z|}:
y,z\in Y,\ \|y\|=\|z\|=1
\right\}.
\]
A possible starting point is its orthogonality-constrained
version
\[
D_{\mathrm{BI}}(Y,|\cdot|)
=
\sup\left\{
\frac{|y|}{|z|}:
\begin{array}{l}
	y,z\in Y,\quad \|y\|=\|z\|=1,\\
	y\mut{\|\cdot\|}z,\quad
	y\iso{\|\cdot\|}z
\end{array}
\right\}.
\]

\section*{Acknowledgements}
Thanks to all the members of the Functional Analysis Research team of the College of Mathematics and Statistics of Anqing Normal University for their discussion and correction of the difficulties and errors encountered in this paper.

\section*{Declarations}
\noindent\textbf{Competing interests.} The authors declare that there are no competing interests.

\medskip\noindent\textbf{Funding.} This research received no external funding.

\medskip\noindent\textbf{Author contributions.} The authors contributed equally to this work.

\medskip\noindent\textbf{Data availability.} No datasets were generated or analysed during this study.

\end{document}